\documentclass[11pt,reqno]{amsart}
\usepackage[T1]{fontenc}
\usepackage{lmodern}
\usepackage{xcolor}
\newif\ifrevisionblue
\revisionbluetrue

\usepackage{amsmath,amssymb,amsthm,mathtools}
\usepackage[protrusion=true,expansion=false]{microtype}
\usepackage[margin=1.15in]{geometry}
\usepackage[colorlinks=true,linkcolor=blue,citecolor=blue,urlcolor=blue]{hyperref}
\allowdisplaybreaks[1]
\numberwithin{equation}{section}

\newtheorem{theorem}{Theorem}[section]
\newtheorem{proposition}[theorem]{Proposition}
\newtheorem{lemma}[theorem]{Lemma}
\newtheorem{corollary}[theorem]{Corollary}
\theoremstyle{remark}
\newtheorem{remark}[theorem]{Remark}
\newcommand{\R}{\mathbb R}
\newcommand{\Id}{\mathrm I}
\newcommand{\dx}{\,\mathrm dx}
\newcommand{\ds}{\,\mathrm ds}
\newcommand{\diver}{\operatorname{div}}

\title[Stationary domains for Robin torsion]{Rigidity of Stationary Domains for Robin Torsion\\ in the Plane}

\author{Qinfeng Li}
\address{School of Mathematics, Hunan University, Changsha, P.R. China.}
\email{liqinfeng1989@gmail.com}
\author{Hang Yang}
\address{School of Mathematics, Hunan University, Changsha, P.R. China.}
\email{hangyang0925@gmail.com}

\subjclass[2020]{35N25, 35B50, 49Q10, 35J25}
\keywords{Robin torsional rigidity, stationary domains,
overdetermined problems, shape optimization,
Airy stress function}

\begin{document}

\begin{abstract}
For every positive Robin parameter, we prove that a bounded simply connected
planar domain with $C^{2,\alpha}$ boundary is stationary for Robin torsional
rigidity under an area constraint if and only if it is a disk. No convexity,
curvature bound, or additional condition on the overdetermined constant is
assumed.  The difficulty is that the stationarity condition couples the boundary trace of the Robin torsion function, its tangential derivative, and the boundary curvature. We encode this condition through a normalized Airy stress function and construct a superharmonic auxiliary function by a quartic minimization calibrated on disks. The boundary structure of this auxiliary function allows us to identify a point where it attains its minimum on the boundary and has a nonnegative outward normal derivative. The strong maximum principle and Hopf's boundary point lemma then force the auxiliary function to be constant. Consequently, for any stationary domain in the class under consideration, the Robin torsion function has a constant boundary trace, and the stationarity condition forces the boundary curvature to be constant.
\end{abstract}

\maketitle

\section{Introduction}
\label{sec:introduction}
Let $\Omega\subset\R^n$, $n\ge2$, be a bounded connected domain with Lipschitz boundary, and let $\beta>0$. The Robin torsion function
$u_{\Omega,\beta}$ is the unique weak solution of
\begin{equation}\label{eq:state}
\begin{cases}
 -\Delta u=1 & \text{in }\Omega,\\
 u_\nu+\beta u=0 & \text{on }\partial\Omega,
\end{cases}
\end{equation}
where $\nu$ is the outward unit normal and
$u_\nu=\partial_\nu u=\nabla u\cdot\nu$ is the normal derivative.
We write $u=u_{\Omega,\beta}$ when the domain and parameter are fixed.
The Robin torsional rigidity is
\[
 T_\beta(\Omega):=\int_\Omega u_{\Omega,\beta}\,\dx.
\]
It represents the total steady-state temperature for a uniform unit heat
source with convective heat transfer through the boundary.
We use $|\Omega|$ for volume (area in the plane), $|\partial\Omega|$ for
boundary measure, and $\ds$ for integration against that measure.

For $\beta>0$, the Saint-Venant inequality,
\begin{align}\label{eq:SV}
T_{\beta}(\Omega)\leq T_{\beta}(\Omega^{\sharp}) \qquad \Omega^{\sharp}\text{ is a ball of volume $|\Omega|$}, 
\end{align}
asserts that
the ball maximizes Robin torsional rigidity among bounded Lipschitz
domains of the same volume. Bucur and Giacomini first established this inequality
through a free discontinuity approach~\cite{BucurGiacomini2015a}; see also the 
survey~\cite{BucurGiacomini2015b}. 
Alvino, Nitsch, and Trombetti~\cite{AlvinoNitschTrombetti2023}
established Talenti-type comparison results for elliptic
problems with Robin boundary conditions, which in particular
yield \eqref{eq:SV}.
Masiello and Paoli~\cite{MasielloPaoli2023} subsequently
characterized the equality case in the planar Talenti-type
comparison. Quantitative stability of \eqref{eq:SV} was established
 in \cite{BucurGiacominiNahon2021, AmatoBarbatoCitoMasielloPaoli2026}. Related optimization problems over the class of tangential polygons have also been studied; see \cite{GuiHuLiYang2026} and the references therein. For a broader overview of Robin problems and related shape
optimization questions, we refer to the survey
chapter by Bucur, Freitas, and Kennedy~
\cite{BucurFreitasKennedy2017}.

Although the Saint-Venant inequality \eqref{eq:SV} identifies balls as global maximizers of Robin torsional rigidity under a volume constraint, it does not rule out the existence of other stationary domains. Throughout, we keep \(\beta>0\) fixed and say that \(\Omega\) is stationary if
\[
\left.\frac{d}{dt}\right|_{t=0}
T_\beta((\Id+tV)(\Omega))=0
\]
for every \(V\in C_c^\infty(\R^n;\R^n)\) satisfying
\(\int_{\partial\Omega}V\cdot\nu\,\ds=0\).
Here \(\Id(x)=x\), \(t\) is the deformation parameter, and the integral condition expresses volume preservation to first order; see \cite[Chapter~5]{HenrotPierre2018}.

For \(C^{2,\alpha}\) domains, \(0<\alpha<1\), the first-variation formula~\cite{BandleWagner2015,LiYang2023} shows that stationarity is equivalent to the existence of a constant \(C_\beta(\Omega)\in\R\) such that the Robin torsion function \(u\), in addition to satisfying \eqref{eq:state}, fulfills
\begin{equation}\label{eq:robin-stationary-n}
-\beta^2u^2+\frac12|\nabla u|^2
+\frac{\beta}{2}Hu^2-u=C_\beta(\Omega)
\qquad\text{on }\partial\Omega.
\end{equation}
Here \(H\) denotes the sum of the principal curvatures with respect to the outward normal, normalized so that \(H=(n-1)/R\) on a sphere of radius \(R\). The constant \(C_\beta(\Omega)\) arises as a Lagrange multiplier for the volume constraint and depends on \(\Omega\) and \(\beta\); its value is not prescribed a priori.
Every ball is stationary, but whether the converse holds is the stationary-domain rigidity question. This problem was raised by Bandle and Wagner~\cite[Example~4.1]{BandleWagner2015} in the context of Robin shape variations. In the heat-transfer setting, Li and Yang~\cite[Corollary~3.2]{LiYang2023} derived the corresponding overdetermined system and explicitly raised the rigidity question for \(f\equiv1\) in \cite[Remark~3.3]{LiYang2023}. This question is the starting point of the present paper; see also \cite[Conjecture~1.1]{LiWeiYao2025} and \cite[Section~1.1]{BarbatoMolinarolo2026}.

To place the Robin stationary-domain question in context, we recall the classical Dirichlet counterpart, Serrin's overdetermined torsion problem:
\begin{equation}\label{eq:serrin}
\begin{cases}
-\Delta u=1 & \text{in }\Omega,\\
u=0 & \text{on }\partial\Omega,\\
u_\nu=\text{constant} & \text{on }\partial\Omega.
\end{cases}
\end{equation}
In the seminal paper~\cite{Serrin1971}, Serrin proved by the method of moving planes that any bounded $C^2$ domain admitting a classical solution to \eqref{eq:serrin} must be a ball. Weinberger~\cite{Weinberger1971} subsequently gave a simpler proof based on an auxiliary $P$-function and the Pohozaev identity. Other approaches include a duality formulation based on mean-value identities for harmonic functions~\cite{PayneSchaefer1989}, shape-derivative arguments combined with Alexandrov's constant-mean-curvature theorem~\cite{ChoulliHenrot1998}, and integral arguments combining Pohozaev-type identities with Newton's inequalities and their rigid equality cases~\cite{BNST08}. For a comprehensive survey of these methods and related developments, see~\cite{NitschTrombetti2018} and the references therein.
Serrin's theorem has motivated extensive work on elliptic problems; see, for instance~\cite{AAB99,ABM24,BianchiniCiraolo18,BNST08,FH02,CS09,CV19,DMAD23,FS15,FM14,FK08,FR22,Fragala12,FFB06,FRS24,HG98,KP98,LXY26,MG20,MP20} and the references therein.
Recent advances have extended Serrin's rigidity theorem to nonsmooth domains under suitable weak formulations. Figalli and Zhang~\cite{FigalliZhang2025}
established rigidity in a finite-perimeter framework under a
uniform upper surface-density bound, including Lipschitz domains,
while Dong and Zhang~\cite{DongZhang2025} provided a different proof
in the Lipschitz setting. Domingo-Pasarin and
Ros-Oton~\cite{DomingoPasarinRosOton2026} obtained an alternative
proof in this setting through boundary regularity for weak
Alt--Caffarelli-type problems. Zhang~\cite{Zhang2026} further showed that, without suitable geometric regularity assumptions, the distributional formulation of Serrin's overdetermined problem admits non-ball counterexamples, while rigidity is recovered in the planar Smirnov class.

Despite these advances in the Dirichlet setting, the Robin stationary-domain problem presents a distinct boundary difficulty. The Robin condition \(u_\nu=-\beta u\) leaves the boundary trace unknown, so the classical moving-plane argument does not directly provide the boundary comparisons available in the Dirichlet setting. Similarly, although Weinberger's \(P\)-function remains subharmonic, its boundary values are no longer constant, and controlling its normal derivative at a boundary maximum generally requires additional information not directly provided by the natural stationarity condition. Existing integral-identity approaches establish rigidity under additional geometric and scalar assumptions~\cite{GavitoneMolinarolo2026,BarbatoMolinarolo2026}. The challenge is therefore to obtain the rigidity from stationarity alone, without imposing such additional conditions.

We now turn to the planar setting of our main rigidity result. Let \(\Omega\subset\mathbb R^2\) be a bounded simply connected domain with \(C^{2,\alpha}\) boundary, where \(0<\alpha<1\). We orient \(\partial\Omega\) clockwise and parametrize it by arclength \(s\). Denoting the unit tangent and outward unit normal by \(\tau\) and \(\nu\), respectively, we define the signed curvature \(\kappa\) through the Frenet formulas
\begin{equation}\label{eq:frenet}
\partial_s\tau=-\kappa\nu,
\qquad
\partial_s\nu=\kappa\tau.
\end{equation}
With this convention, \(\kappa\geq0\) for convex domains and \(\kappa=1/R\) on the boundary of a disk of radius \(R\). For a sufficiently smooth function \(f\), we write \(f_s=\partial_s(f|_{\partial\Omega})\) for its tangential derivative. Since \(H=\kappa\) in two dimensions, the overdetermined system~\eqref{eq:robin-stationary-n} takes the form
\begin{align}
\label{eq:system}
\begin{cases}
-\Delta u=1 & \text{in } \Omega,\\
\partial_\nu u+\beta u=0 & \text{on } \partial\Omega,\\
-\beta^2u^2+\dfrac12|\nabla u|^2
+\dfrac{\beta}{2}\kappa u^2-u=C_{\beta}(\Omega)
& \text{on } \partial\Omega.
\end{cases}
\end{align}
Here $C_\beta(\Omega)$ is the same a priori unknown constant associated
with the volume constraint; no value or comparison with a model value is
assumed.

Our main result establishes rigidity of \eqref{eq:system} for every positive Robin parameter.
\begin{theorem}\label{thm:main}
Let $0<\alpha<1$ and $\beta>0$, and let $\Omega\subset\R^2$ be a
nonempty bounded simply connected domain with $C^{2,\alpha}$ boundary. Then \eqref{eq:system} admits a solution $u\in C^{2,\alpha}(\overline{\Omega})$ if and only if $\Omega$ is a disk.
\end{theorem}
%=======
The additional assumptions in two recent higher-dimensional Robin rigidity results can be expressed in terms of
\[
R_\Omega:=\frac{n|\Omega|}{P(\Omega)},\qquad
C_\beta^*(\Omega):=-\frac{1}{2n^2}
\left(R_\Omega^2+\frac{(n+1)R_\Omega}{\beta}\right),
\qquad
\kappa_{\min}:=\min_{\substack{x\in\partial\Omega\\1\le i\le n-1}}\kappa_i(x),
\]
where \(P(\Omega)=|\partial\Omega|\), \(\kappa_i\) are the principal curvatures of \(\partial\Omega\), and \(C_\beta^*(\Omega)\) is the stationary constant for a ball of radius \(R_\Omega\) in our normalization.
More precisely, \cite[Theorem~1]{GavitoneMolinarolo2026} establishes rigidity under
\(\kappa_{\min}+\beta\ge0\) and \(C_\beta(\Omega)=C_\beta^*(\Omega)\), whereas the \(p=2\) case of \cite[Theorem~1.1]{BarbatoMolinarolo2026} assumes
\(\kappa_{\min}+\beta>0\) and \(C_\beta(\Omega)\ge C_\beta^*(\Omega)\).
The role of these assumptions is clarified by the \(p=2\) specialization of the integral identity in \cite[Theorem~4.4]{BarbatoMolinarolo2026}, which, in our normalization, reads
\small{
\[
\begin{aligned}
\int_\Omega\left|D^2u+\frac1n\Id\right|^2\,dx
+\beta\int_{\partial\Omega}
\left(u_\nu+\frac{R_\Omega}{n}\right)^2\,dS=2\beta P(\Omega)
\bigl(C_\beta^*(\Omega)-C_\beta(\Omega)\bigr)-\int_{\partial\Omega}
\left[\mathrm{II}(\nabla_\tau u,\nabla_\tau u)
+\beta|\nabla_\tau u|^2\right]\,dS,
\end{aligned}
\]}
where \(\mathrm{II}\) is the second fundamental form and \(\nabla_\tau u\) is the tangential gradient along \(\partial\Omega\). 
Since
\[
\mathrm{II}(\nabla_\tau u,\nabla_\tau u)
+\beta|\nabla_\tau u|^2
\ge(\kappa_{\min}+\beta)|\nabla_\tau u|^2, 
\]
the conditions \(\kappa_{\min}+\beta\ge0\) and \(C_\beta(\Omega)\ge C_\beta^*(\Omega)\) force both nonnegative terms on the left-hand side to vanish. In particular, \(D^2u=-\frac1n\Id\), and the Robin boundary condition then implies that \(\Omega\) is a ball.
Thus, these integral-identity approaches obtain rigidity by combining stationarity with additional curvature and scalar constraints. 

In contrast, the stationary-domain conjecture
requires only that the boundary expression in
\eqref{eq:system} be constant; it neither
prescribes the value of this constant nor
imposes any curvature bound. In this paper,
Theorem~\ref{thm:main} resolves the
stationary-domain rigidity problem for all
\(\beta>0\) among simply connected planar
domains, without either a curvature bound or
an additional assumption on \(C_\beta(\Omega)\).

The proof of Theorem~\ref{thm:main} is inspired by rigidity methods that exploit differential identities and tensorial structures to characterize equality cases~\cite{Obata1971,MaWu2024,GuiLiWangXu2026}.  We seek a differential identity whose equality case characterizes the torsion profile of a disk. We implement this idea through a normalized Airy stress function and a quartic envelope calibrated on disks, obtaining a  rigidity mechanism distinct from the integral-identity criteria above. The precise connections are explained in Remark~\ref{rem:tensor-comparison}.

The proof rests on three parts.
First, we use the divergence-free stress tensor~\eqref{eq:stress} to express the stationarity condition as the boundary flux identity~\eqref{eq:T-nu}. In the planar simply connected setting, the associated normalized Airy stress function \(W\) can be globally defined and normalized, using the connectedness of \(\partial\Omega\), so that
\[
\Delta W=k:=u+C_\beta(\Omega)\quad\text{in }\Omega,
\qquad
W=0,\quad \nabla W=\frac{\beta}{4}u^2\nu
\quad\text{on }\partial\Omega.
\]
The advantage of this normalization is that the boundary data involve only \(u\), with no explicit dependence on its tangential derivative or the boundary curvature (Lemma~\ref{lem:airy}).

Second, we construct an auxiliary function $b$ which is constant on the stationary domain whenever it is a disk. 
The identity
\[
\nabla W=-\left(k+\frac12|\nabla u|^2\right)\nabla u.
\]
valid on disks, motivates the auxiliary function
\[
b=W+\min_{v\in\R^2}
\left\{
\frac14|v|^4+k|v|^2-2v\cdot\nabla W
\right\}.
\]
On $\partial\Omega$, the auxiliary function admits the
representation $b=\mathcal B(u)$, where $\mathcal B$ is
defined in \eqref{eq:B} and is strictly concave on $(0,\infty)$.
Together with the stationarity condition and the disk calibration, this shows that $b$ is constant
if and only if $\Omega$ is a disk
(Proposition~\ref{prop:constant-envelope}).
It therefore remains to establish the constancy of $b$.

Finally, we establish the superharmonicity of $b$
by combining a differential identity with a viscosity
argument at the critical points of $W$
(Proposition~\ref{prop:source} and
Lemma~\ref{lem:critical}).
The strict concavity of $\mathcal B$ ensures that
a boundary minimum of $b$ is attained at a point
$P$ where $u$ achieves its boundary minimum or maximum.
Hence $u_s(P)=0$, and the normal derivative formula
gives $b_\nu(P)\ge0$ (Lemma~\ref{lem:boundary}).
If $b$ were nonconstant, the strong minimum principle
and Hopf's boundary point lemma would instead yield
$b_\nu(P)<0$, a contradiction.
Thus $b$ is constant, and hence $\Omega$ is a disk.
As a complementary result, the equality case of the
differential identity also yields
$D^2u=-\frac12\Id$
(Remark~\ref{rem:profile-rigidity}).

The method also extends to semilinear Robin problems
with positive, nonincreasing source terms. More precisely,
we establish a corresponding rigidity result for
$-\Delta u=g(u)$ with $g>0$ and $g'\le0$
(Theorem~\ref{thm:semilinear}), addressing the planar
case of the question raised in
\cite[Section~4.1, Example~2]{BandleWagner2015}.

\textbf{Organization.}
Section~\ref{sec:preliminaries} constructs the normalized potential.
Sections~\ref{sec:scalar} and~\ref{sec:critical} introduce the envelope and
prove its superharmonicity. Section~\ref{sec:boundary} completes the
boundary argument. Section~\ref{sec:semilinear} extends the method
to semilinear equations with positive, nonincreasing source terms. The complementary $P$-function proof in
Section~\ref{sec:P-function} permits multiply connected domains under
a curvature condition, without a scalar restriction on $C_\beta(\Omega)$.

\section{Stress tensor and the normalized Airy stress function}
\label{sec:preliminaries}
We begin by recalling some basic properties of the Robin torsion function. For \(\beta>0\) and a bounded \(C^{2,\alpha}\) domain \(\Omega\subset\R^2\), standard elliptic theory and Schauder estimates ensure that \eqref{eq:state} admits a unique solution \(u_{\Omega,\beta}\in C^{2,\alpha}(\overline\Omega)\); see \cite[Section~6.7, Theorem~6.31]{GilbargTrudinger2001}. The maximum principle and Hopf's boundary point lemma yield \(u_{\Omega,\beta}>0\) on \(\overline\Omega\). In what follows, we simply write \(u=u_{\Omega,\beta}\) whenever no confusion arises.

The weak formulation of \eqref{eq:state} is
\[
\int_\Omega\nabla u\cdot\nabla\varphi\,\dx
+\beta\int_{\partial\Omega}u\varphi\,\ds
=\int_\Omega\varphi\,\dx,
\qquad \forall\,\varphi\in H^1(\Omega).
\]
In particular, choosing \(\varphi\equiv1\) and \(\varphi=u\) gives
\begin{equation}\label{eq:mass-energy}
|\Omega|=\beta\int_{\partial\Omega}u\,\ds,
\qquad
\int_\Omega u\,\dx
=\int_\Omega|\nabla u|^2\,\dx
+\beta\int_{\partial\Omega}u^2\,\ds.
\end{equation}

For the tensor calculations below, we adopt the following conventions. All vector norms are Euclidean, while the norm of a \(2\times2\) matrix \(A\) is the Frobenius norm \(|A|^2=\sum_{i,j}A_{ij}^2\). For \(a,b,\xi,\eta\in\R^2\), we write \((a\otimes b)_{ij}=a_i b_j\) and \(A(\xi,\eta)=\xi\cdot A\eta\). The divergence of a matrix field is taken row-wise, namely \((\diver A)_i=\sum_j\partial_jA_{ij}\).

To determine \(C_\beta(\Omega)\) and construct the normalized Airy stress function, we introduce the stress-energy tensor associated with the Lagrangian \(L(u,\nabla u)=|\nabla u|^2/2-u\):
\begin{equation}\label{eq:stress}
\mathsf T=\left(\frac12|\nabla u|^2-u\right)\Id
-\nabla u\otimes\nabla u.
\end{equation}
Such tensors are standard in the variational theory of elliptic equations and have been used to derive conservation laws, monotonicity formulas, and Pohozaev-type identities; see, for instance, \cite{PucciSerrin1986,Alikakos2011,AlikakosFaliagas2012}. For the torsion equation, \(\operatorname{div}\mathsf T=0\) follows directly by differentiation. The boundary flux \(\mathsf T\nu\) is analogous to the traction vector in continuum mechanics~\cite{Gurtin1981}. The Robin condition and stationarity further imply that the shifted flux \((\mathsf T-C_\beta(\Omega)\Id)\nu\) is an arclength derivative, providing the boundary structure needed to normalize the associated Airy stress function.

\begin{lemma}\label{lem:multiplier}
Let $\Omega\subset\mathbb R^2$ be a bounded connected domain
with $C^{2,\alpha}$ boundary, where $\alpha\in(0,1)$, and let
$\beta>0$. Suppose that $u\in C^{2,\alpha}(\overline\Omega)$ is the solution to \eqref{eq:system}, and let the tensor $\mathsf T$ be defined by \eqref{eq:stress}. Denote the constant in \eqref{eq:system} by $C_{\beta}(\Omega)$. Then
\begin{equation}\label{eq:stress-conservation}
 \diver\mathsf T=0,\qquad \operatorname{tr}\mathsf T=-2u \qquad \text{in }\Omega,
\end{equation}
and
\begin{equation}\label{eq:T-nu}
 \mathsf T\nu =C_{\beta}(\Omega)\nu+\frac{\beta}{2}(u^{2}\tau)_{s} \qquad \text{on }\partial \Omega,
\end{equation}
where $\tau$ and $\nu$ are the unit tangent and outward unit  normal to $\partial \Omega$. Moreover, we have
\begin{equation}\label{eq:constant}
C_{\beta}(\Omega)=-\frac{1}{|\Omega|}\Big(\int_\Omega|\nabla u|^2\dx +\frac{3\beta}{4}\int_{\partial \Omega} u^2\ds\Big)<0.
\end{equation}
\end{lemma}

\begin{proof}
Since $\Delta u=-1$, direct differentiation gives
\[
\diver\mathsf T=-(1+\Delta u)\nabla u=0 \qquad \text{in }\Omega.
\]
Taking the trace yields
$\operatorname{tr}\mathsf T=-2u$. Hence
\(
\diver(\mathsf T x)=-2u
\) in $\Omega$.
On the boundary, the Robin condition and the last equation in \eqref{eq:system} become
\begin{equation}\label{eq:scalar}
\nabla u=u_{s}\tau-\beta u\nu,
\qquad
\frac12 u_{s}^2+\frac\beta2(\kappa-\beta)u^2-u=C_{\beta}(\Omega).
\end{equation}
Consequently,
\begin{align*}
\mathsf T\nu
&=
\left(\frac12 u_s^2-\frac{\beta^2}{2}u^2-u\right)\nu
+\frac{\beta}{2}(u^2)_s\tau\\
&=\left(C_{\beta}(\Omega)-\frac\beta2\kappa u^2\right)\nu+\frac\beta2(u^2)_s\tau\\
&=C_{\beta}(\Omega)\nu+\frac{\beta}{2}(u^{2}\tau)_{s}.
\end{align*}

Using \eqref{eq:T-nu} and the divergence theorem, we obtain
\begin{align}\label{eq:C1}
 -2\int_\Omega u\dx=\int_{\Omega} \diver (\mathsf T x)dx=\int_{\partial \Omega}(\mathsf T\nu)\cdot x\ds=2|\Omega|C_{\beta}(\Omega)+\frac{\beta}{2}\int_{\partial \Omega}(u^{2}\tau)_{s} \cdot x ds,
\end{align}
where we used $(\mathsf T\nu)\cdot x= (\mathsf T x)\cdot \nu$ since $\mathsf T$ is symmetric.

Since $\partial \Omega$ is closed and $x_{s}=\tau$, integration by parts yields
\begin{align}\label{eq:C2}
\int_{\partial \Omega}(u^{2}\tau)_{s}\cdot x ds=-\int_{\partial \Omega}u^{2}\tau\cdot\tau ds=-\int_{\partial \Omega}u^{2} ds.
\end{align}
Combining \eqref{eq:mass-energy}, \eqref{eq:C1}, and \eqref{eq:C2} proves
\eqref{eq:constant}. 
\end{proof}

The classical Airy stress-function
representation, originating in planar elasticity~\cite{Airy1863},
allows a symmetric divergence-free tensor on a simply connected
planar domain to be expressed in terms of the second derivatives
of a scalar function; see
\cite[Theorem~2.2]{CiarletGeymonatKrasucki2013} for a rigorous
formulation, and \cite{GeymonatKrasucki2000} for Lipschitz
domains. 
By \eqref{eq:stress-conservation}, the tensor
$\mathsf S:=\mathsf T-C_\beta(\Omega)\Id$ is symmetric and
divergence-free. Since $\Omega$ is simply connected, it admits
an \emph{Airy stress function} $\phi$ such that, in Cartesian
coordinates,
\[
\mathsf S=
\begin{pmatrix}
\phi_{yy}&-\phi_{xy}\\
-\phi_{xy}&\phi_{xx}
\end{pmatrix}.
\]
The function $\phi$ is unique up to an affine term
$\gamma+a\cdot x$, with $\gamma\in\R$ and $a\in\R^2$,
which leaves its Hessian unchanged.
Let $J$ denote the clockwise quarter-turn rotation,
so that $\tau=J\nu$. We introduce $W=-\phi/2$, for which
\[
D^2W=-\frac12J\mathsf S J^{\mathsf T},
\]
which gives \eqref{eq:airy-hessian}. Moreover,
\eqref{eq:T-nu} implies
\[
\partial_s\left(\nabla W-\frac{\beta}{4}u^2\nu\right)=0
\qquad\text{on }\partial\Omega.
\]
Since $\partial\Omega$ is connected, the affine freedom
of $W$ allows us to impose the boundary normalization
in \eqref{eq:airy-boundary}. We refer to the resulting
uniquely determined function as the normalized
Airy stress function. The precise statement is given in
Lemma~\ref{lem:airy}.

\begin{lemma}\label{lem:airy}
Let $\Omega\subset\mathbb R^2$ be a bounded simply
connected domain with $C^{2,\alpha}$ boundary,
where $0<\alpha<1$, and let $\beta>0$. Suppose that $u\in C^{2,\alpha}(\overline\Omega)$ is the solution to \eqref{eq:system}. There is a unique $W\in C^{3,\alpha}(\overline\Omega)$ satisfying
\begin{equation}\label{eq:airy-hessian}
 D^2W=\frac12\left[
 \left(k+\frac12|\nabla u|^2\right)\Id
 -\nabla u\otimes\nabla u\right],
\end{equation}
where $k(x)=u(x)+C_{\beta}(\Omega)$, and the normalization
\begin{equation}\label{eq:airy-boundary}
\begin{cases}
 \Delta W=k \quad &\text{in }\Omega,\\
 W=0 &\text{on }\partial \Omega,\\
 W_\nu=\frac\beta4 u^2 &\text{on }\partial \Omega.
\end{cases}
\end{equation}
\end{lemma}

\begin{proof}

Let $M$ denote the matrix on the right-hand side of
\eqref{eq:airy-hessian}. Since
$u\in C^{2,\alpha}(\overline\Omega)$, we have
$M\in C^{1,\alpha}(\overline\Omega)$. Direct differentiation,
using $\Delta u=-1$, gives
\[
 \partial_yM_{11}-\partial_xM_{12}
 =\frac{u_y}{2}(1+\Delta u)=0,
 \qquad
 \partial_yM_{21}-\partial_xM_{22}
 =-\frac{u_x}{2}(1+\Delta u)=0.
\]
Since \(\Omega\) is simply connected, the vanishing curl identities imply path independence of the corresponding line integrals. Hence there exist
$V_1,V_2\in C^2(\Omega)$ such that
\[
 \partial_jV_i=M_{ij},\qquad i,j=1,2.
\]
By the symmetry of $M$,
\[
 \partial_yV_1-\partial_xV_2=M_{12}-M_{21}=0.
\]
Applying the same criterion once more, we obtain
$W\in C^3(\Omega)$ satisfying
\[
 \nabla W=(V_1,V_2),\qquad D^2W=M.
\]
Thus $W$ is unique up to an affine function.

The regularity extends to $\overline\Omega$. Indeed, in each sufficiently
small Lipschitz boundary chart, interior points $x,y$ can be joined inside
$\Omega$ by a curve of length at most $K|x-y|$. Integrating a continuous
gradient along such curves gives a continuous extension of its primitive;
if the gradient is $C^{0,\alpha}$, subtracting its value at $x$ in the
integral gives a Taylor remainder bounded by $C|x-y|^{1+\alpha}$.
Applying this successively to $\nabla V_i=(M_{i1},M_{i2})$ and
$\nabla W=(V_1,V_2)$ proves $V_i\in C^{2,\alpha}(\overline\Omega)$ and
$W\in C^{3,\alpha}(\overline\Omega)$.

Finally, $W$ is determined up to an affine function,
since two potentials with the same Hessian differ
by an affine function on the connected domain $\Omega$.
On $\partial \Omega$, the decomposition $\nabla u=u_s\tau-\beta u\nu$ and
\eqref{eq:scalar} yield
\[
 M(\tau,\tau)
 =\frac12\left[k+\frac{\beta^2u^2-u_s^2}{2}\right]
 =\frac\beta4\kappa u^2,
 \qquad M(\tau,\nu)=\frac\beta2u u_s.
\]
Using \eqref{eq:frenet}, we obtain
\[
 \partial_s\nabla W=M\tau
 =\frac\beta4\kappa u^2\tau+\frac\beta2u u_s\nu
 =\partial_s\left(\frac\beta4u^2\nu\right).
\]
Because $\partial \Omega$ is connected, the difference of these two vectors is
a single constant vector $a_0$. Replacing $W$ by $W-a_0\cdot x$
makes its boundary gradient equal to $(\beta u^2/4)\nu$.
Its tangential derivative is then zero, so subtracting one more
constant gives $W=0$ on $\partial \Omega$. Taking the trace of \(\eqref{eq:airy-hessian}\) yields \(\Delta W=k\) in \(\Omega\). The boundary gradient identity established above gives \(W_\nu=\beta u^2/4>0\) on \(\partial\Omega\). Finally, the difference of two normalized
potentials is affine, with zero value and zero gradient on $\partial \Omega$,
and hence vanishes identically.
\end{proof}

\begin{remark}
Simple connectivity is used for the global integrations defining $W$;
connectedness of $\partial\Omega$ makes the affine normalization common
to the whole boundary. These are requirements of the present argument;
we do not claim that simple connectivity is necessary for rigidity. 
\end{remark}

\section{Disk calibration and the quartic envelope}
\label{sec:scalar}

Let $\Omega=B_R\subset\R^2$ be the disk of radius
$R>0$ centered at the origin, and set $\rho=|x|$.
The Robin torsion function and the boundary constant are
\begin{equation}\label{eq:disk}
u(x)=\frac{R^2-\rho^2}{4}+\frac{R}{2\beta},
\qquad
C_{\beta}(B_R)=-\frac{R^2}{8}-\frac{3R}{8\beta}.
\end{equation}
The corresponding normalized Airy stress function is
\begin{equation}\label{eq:disk-airy}
W(x)=\frac{(\rho^2-R^2)(R^2+2R/\beta-\rho^2)}{64}.
\end{equation}
Writing $k=u+C_{\beta}(B_R)$, we have
\[
|\nabla u|^2(x)=\frac{\rho^2}{4},
\qquad
k(x)=\frac{R^2+R/\beta}{8}-\frac{\rho^2}{4}.
\]
A direct calculation then gives
\begin{equation}\label{eq:disk-cancellation}
W-k|\nabla u|^2-\frac34|\nabla u|^4
=-\frac{R^2(R^2+2R/\beta)}{64}
\qquad\text{in }B_R.
\end{equation}

To extend this cancellation to general stationary domains
in a form suited to the Laplacian estimates, we
express the correction as a pointwise minimum over an
independent vector. The key observation is that
\begin{equation}\label{eq:disk-gradient-calibration}
\nabla W=-\left(k+\frac12|\nabla u|^2\right)\nabla u \qquad\text{in }B_{R}.
\end{equation}
Set $v_0=-\nabla u$ and
\[
A:=k+\frac12|\nabla u|^2
=\frac{R^2+R/\beta-\rho^2}{8}>0.
\]
Then $\nabla W=Av_0$. For each fixed $x$, consider
\[
F_x(v)=\frac14|v|^4+k(x)|v|^2
       -2v\cdot\nabla W(x),
\qquad v\in\R^2.
\]
Its critical-point equation is
$\nabla W=(k+|v|^2/2)v$.
Since $A=k+|v_0|^2/2$, a direct expansion gives
\[
F_x(v)-F_x(v_0)
=A|v-v_0|^2
+\frac14\bigl(|v|^2-|v_0|^2\bigr)^2\ge0.
\]
Thus $v_0$ is the unique global minimizer, and
\[
\min_{v\in\R^2}F_x(v)
=F_x(v_0)
=-k|\nabla u|^2-\frac34|\nabla u|^4.
\]

Motivated by this representation, for any bounded simply
connected $C^{2,\alpha}$ domain $\Omega$ admitting a solution
of \eqref{eq:system}, we define
\begin{equation}\label{eq:envelope}
b(x)=W(x)+\min_{v\in\R^2}
\left\{\frac14|v|^4+k(x)|v|^2
-2v\cdot\nabla W(x)\right\},
\qquad x\in\overline\Omega,
\end{equation}
where $W$ and $k$ are given by Lemma~\ref{lem:airy}.
The next proposition shows that the
constancy of $b$ characterizes disks.

\begin{proposition}\label{prop:constant-envelope}
Under the assumptions of Lemma~\ref{lem:airy}, let $b$ be
defined by \eqref{eq:envelope}, and put $C=C_\beta(\Omega)$.
For $r>0$, define
\begin{equation}\label{eq:B}
\mathcal B(r):=\min_{s\ge0}
\left\{\frac14s^4+(r+C)s^2-\frac{\beta}{2}r^2s\right\}.
\end{equation}
Then $\mathcal B$ is strictly concave on $(0,\infty)$ and
$b=\mathcal B(u)$ on $\partial\Omega$.
Moreover, $b$ is constant on $\overline{\Omega}$ if and only if $\Omega$ is a disk.
\end{proposition}

\begin{proof}

On $\partial\Omega$, the normalization \eqref{eq:airy-boundary}
gives $W=0$ and $\nabla W=(\beta/4)u^2\nu$.
For each fixed $|v|=s$, the expression minimized in
\eqref{eq:envelope} is smallest when $v=s\nu$.
This yields
\begin{equation}\label{eq:boundary-profile}
b=\mathcal B(u)\qquad\text{on }\partial\Omega,
\end{equation}
where $\mathcal B$ is defined in \eqref{eq:B}.
Set
\[
F(r,s)=\frac14s^4+(r+C)s^2-\frac{\beta}{2}r^2s,
\qquad r>0,\ s\ge0.
\]
For each $r>0$, $F(r,s)\to+\infty$ as $s\to\infty$,
while $\partial_sF(r,0)=-\beta r^2/2<0$.
Thus $F(r,\cdot)$ attains its minimum at some $s_r>0$.
Let $r_1\ne r_2$ be positive, $0<\theta<1$, and
$r_\theta=\theta r_1+(1-\theta)r_2$.
Choosing a minimizer $s_\theta=s_{r_\theta}>0$, we obtain
\[
\begin{aligned}
\mathcal B(r_\theta)
&=\theta F(r_1,s_\theta)
 +(1-\theta)F(r_2,s_\theta)
 +\frac{\beta}{2}s_\theta\theta(1-\theta)(r_1-r_2)^2\\
&\ge\theta\mathcal B(r_1)+(1-\theta)\mathcal B(r_2)
 +\frac{\beta}{2}s_\theta\theta(1-\theta)(r_1-r_2)^2\\
&>\theta\mathcal B(r_1)+(1-\theta)\mathcal B(r_2).
\end{aligned}
\]
Hence $\mathcal B$ is strictly concave on $(0,\infty)$.

Let $m=\min_{\partial \Omega}u$, and let $M=\max_{\partial \Omega}u$. Since $\partial\Omega$ is connected,
$u(\partial\Omega)=[m,M]$.
Suppose that $b$ is constant on $\overline\Omega$.
Then $\mathcal B$ is constant on $[m,M]$, and its strict
concavity forces $m=M$. Thus $u$ is constant on
$\partial\Omega$. By the stationarity condition
\eqref{eq:scalar}, the curvature of $\partial\Omega$
is also constant. Since $\partial\Omega$ is a smooth
simple closed curve, it must be a circle, and hence
$\Omega$ is a disk.
Conversely, if $\Omega=B_R$, the disk calibration
preceding \eqref{eq:envelope}, together with
\eqref{eq:disk-cancellation}, shows that $b$ is constant.
\end{proof}

\begin{remark}
The pointwise minimization in \eqref{eq:envelope}
can be expressed in terms of the Legendre--Fenchel
transform; see
\cite[Section~11.A, Eq.~11(1)]{RockafellarWets1998}.
For fixed $k\in\R$, let
$H_k(v)=|v|^4/4+k|v|^2$, with Fenchel conjugate
\[
H_k^*(p)=\sup_{v\in\R^2}
\{p\cdot v-H_k(v)\}.
\]
Then $b=W-H_k^*(2\nabla W)$, where $k=k(x)$. 
\end{remark}

In view of the rigidity characterization Proposition~\ref{prop:constant-envelope}, it remains to show that the existence of a solution to the overdetermined problem forces $b$ to be constant. To this end, we investigate the differential and boundary properties of $b$.

\begin{lemma}\label{lem:envelope}
Under the hypotheses of Lemma~\ref{lem:airy}, let
$u$, $k=u+C_\beta(\Omega)$, and $W$ be as defined therein.
Let $b$ be defined by \eqref{eq:envelope}, and set
$t=|\nabla W|^2$. Then $b\in C(\overline\Omega)$,
and the following properties hold.

\begin{itemize}
\item[(1)] At each point where $t>0$, there is a unique
$\ell=\ell(k,t)>\max\{k,0\}$ satisfying
\begin{equation}\label{eq:root}
t=2\ell^2(\ell-k)=\ell^2\sigma,
\qquad \sigma:=\sigma(k,t)=2(\ell(k,t)-k)>0.
\end{equation}
The minimum in \eqref{eq:envelope} is attained uniquely at
\begin{equation}\label{eq:minimizer}
v_*=\frac{\nabla W}{\ell},
\qquad |v_*|^2=\sigma.
\end{equation}
Moreover,
\begin{equation}\label{eq:b}
b=W+\Phi(k,t),
\qquad
\Phi(k,t)=-k\sigma(k,t)-\frac34\sigma(k,t)^2.
\end{equation}

\item[(2)] At each point where $t=0$,
\begin{equation}\label{eq:extension}
b=W-(-k)_+^2,
\end{equation}
where $a_+:=\max\{a,0\}$. The unique minimizer is
$v=0$ when $k\ge0$, whereas for $k<0$ the set of
minimizers is
\[
\{v\in\mathbb R^2:|v|^2=-2k\}.
\]

\item[(3)] The functions $\ell(k,t)$ and $\sigma(k,t)$
are smooth on $\mathbb R\times(0,\infty)$.
Consequently, $b\in C^2(\Omega\cap\{t>0\})$, and
its derivatives up to order two extend continuously
to $\partial\Omega\cap\{t>0\}$.
\end{itemize}
\end{lemma}
%========

\begin{proof}
Write $w=\nabla W$ and
$f_{k,w}(v)=|v|^4/4+k|v|^2-2v\cdot w$.
Since $k,w$ are bounded on $\overline\Omega$, uniform coercivity gives
$M>0$ such that all global minimizers lie in $\overline B_M$.
They exist by compactness, and their minimum value $m$ satisfies
\[
 |m(x)-m(y)|\le M^2|k(x)-k(y)|+2M|w(x)-w(y)|.
\]
Thus $b=W+m$ is continuous.

If $w\ne0$, testing a small positive multiple of $w$ gives a negative
value, so every minimizer is nonzero. For fixed $|v|$, the linear term
is minimized only when $v$ is a positive multiple of $w$. The critical-point
equation is therefore
\[
 w=\left(k+\frac12|v_*|^2\right)v_*.
\]
Setting $\sigma=|v_*|^2$ and $\ell=k+\sigma/2$ gives
$\ell>\max\{k,0\}$ and $t=2\ell^2(\ell-k)$.
On this interval $G_k(\ell):=2\ell^2(\ell-k)$ increases from zero to
infinity, since
\[
 G_k'(\ell)=2\ell(3\ell-2k)=2\ell(\ell+\sigma)>0.
\]
Thus there is exactly one such $\ell$, and every minimizer is $w/\ell$.
Substitution gives
$f_{k,w}(v_*)=\sigma^2/4+k\sigma-2\ell\sigma
=-k\sigma-3\sigma^2/4$, proving (1).

If $w=0$, minimizing $|v|^4/4+k|v|^2$ gives the unique minimizer $0$
for $k\ge0$. For $k<0$, completing the square gives
$(|v|^2+2k)^2/4-k^2$, whose minimizers have $|v|^2=-2k$.
This proves (2).

Finally, $\partial_\ell(G_k(\ell)-t)=2\ell(\ell+\sigma)>0$, so the
implicit function theorem and uniqueness give smooth $\ell,\sigma,\Phi$
on $\R\times(0,\infty)$. Since $k,t\in C^{2,\alpha}(\overline\Omega)$,
the formula $b=W+\Phi(k,t)$ proves (3) by the chain rule, also at boundary
points where $t$ stays positive.
\end{proof}

\section{Superharmonicity and minimum principles}
\label{sec:critical}

By Proposition~\ref{prop:constant-envelope},
the proof of Theorem~\ref{thm:main} reduces to
showing that $b$ is constant whenever
\eqref{eq:system} admits a solution.
To this end, we first establish the
superharmonicity of $b$ in $\Omega$.
Away from the critical set of $W$, the function
$b$ is $C^2$ by Lemma~\ref{lem:envelope},
and the key ingredient is the full Hessian
identity \eqref{eq:airy-hessian}.

\begin{proposition}\label{prop:source}

Under the assumptions and notation of Lemma~\ref{lem:airy},
let $b$ be defined by \eqref{eq:envelope}.
At each point where $|\nabla W|>0$, let $\ell,\sigma$
be given by \eqref{eq:root}, and set
\[
v_*:=\frac{\nabla W}{\ell},
\qquad
Z_\pm:=|\nabla u\pm v_*|^2.
\]
Then
\begin{equation}\label{eq:source}
\Delta b
=-\frac{Z_+(2\ell Z_++\sigma Z_-+8\ell\sigma)}
{8\ell(\ell+\sigma)}\le0,
\end{equation}
with equality if and only if $\nabla W=-\ell\nabla u$.

\end{proposition}

\begin{proof}
We retain the notation of Lemma~\ref{lem:envelope}. In $\{t>0\}$, Lemma~\ref{lem:envelope} gives
\[
 b=W+\Phi(k,t),\qquad
 \Phi(k,t)=-k\sigma-\frac34\sigma^2,
 \qquad t=\ell^2\sigma,\quad \ell=k+\frac\sigma2.
\]
Subscripts $k,t$ on $\sigma$ and $\Phi$ denote partial derivatives with
respect to these scalar parameters. Differentiating
$t=\sigma(k+\sigma/2)^2$ with respect to $k$ and $t$, respectively, yields
\[
 \sigma_k=-\frac{2\sigma}{\ell+\sigma},\qquad
 \sigma_t=\frac1{\ell(\ell+\sigma)}.
\]
Since $k+3\sigma/2=\ell+\sigma$, differentiation of $\Phi$
now gives
\begin{equation}\label{eq:phi-derivatives}
 \begin{gathered}
 \Phi_k=\sigma,\qquad \Phi_t=-\frac1\ell,\\
 \Phi_{kk}=-\frac{2\sigma}{\ell+\sigma},\qquad
 \Phi_{kt}=\frac1{\ell(\ell+\sigma)},\qquad
 \Phi_{tt}=\frac1{2\ell^3(\ell+\sigma)}.
 \end{gathered}
\end{equation}
In particular, since $\nabla k=\nabla u$,
\begin{equation}\label{eq:differential-b}
 \nabla b=\nabla W+\sigma\nabla u-\frac{\nabla t}{\ell}.
\end{equation}

Next, the Hessian identity \eqref{eq:airy-hessian} and
$t=|\nabla W|^2$ imply
\begin{align}
 \nabla t
 &=\left(k+\frac12|\nabla u|^2\right)\nabla W
       -(\nabla u\cdot\nabla W)\nabla u,
       \label{eq:t-gradient}\\
 \Delta t
 &=k^2+\frac14|\nabla u|^4+2\nabla u\cdot\nabla W,
       \label{eq:t-laplacian}\\
 \nabla u\cdot\nabla t
 &=\left(k-\frac12|\nabla u|^2\right)
       (\nabla u\cdot\nabla W),\notag\\
 |\nabla t|^2
 &=\left(k+\frac12|\nabla u|^2\right)^2t
       -2k(\nabla u\cdot\nabla W)^2.
       \label{eq:t-contractions}
\end{align}
Indeed, $\nabla t=2D^2W\,\nabla W$ gives the first identity,
and taking its scalar product with $\nabla u$ and its squared
norm gives the last two. For the second, use $\Delta W=k$ and
$\nabla k=\nabla u$ to write
\[
 \Delta t=2|D^2W|^2+2\nabla W\cdot\nabla k,
 \qquad
 2|D^2W|^2=k^2+\frac14|\nabla u|^4,
\]
where the latter equality follows directly from
\eqref{eq:airy-hessian} in dimension two.

Using $\Delta W=k$ and $\Delta k=-1$, the second-order
chain rule and \eqref{eq:phi-derivatives} give
\[
 \begin{aligned}
 \Delta b={}&k-\sigma-\frac{\Delta t}{\ell}
     -\frac{2\sigma|\nabla u|^2}{\ell+\sigma}
     +\frac{2\nabla u\cdot\nabla t}{\ell(\ell+\sigma)}
     +\frac{|\nabla t|^2}{2\ell^3(\ell+\sigma)}.
 \end{aligned}
\]
Substitute the identities above and use
$k=\ell-\sigma/2$, $\nabla W=\ell v_*$, and
$|v_*|^2=\sigma$. Expanding and collecting terms yields
\[
 \begin{aligned}
 -8\ell(\ell+\sigma)\Delta b
 ={}&2\ell\bigl(|\nabla u|^2+\sigma
                   +2\nabla u\cdot v_*\bigr)^2\\
 &+\sigma\bigl[\bigl(|\nabla u|^2+\sigma\bigr)^2
                   -4(\nabla u\cdot v_*)^2\bigr]\\
 &+8\ell\sigma\bigl(|\nabla u|^2+\sigma
                   +2\nabla u\cdot v_*\bigr).
 \end{aligned}
\]
By the definitions of $Z_\pm$,
\[
 Z_+=|\nabla u|^2+\sigma+2\nabla u\cdot v_*,\qquad
 Z_+Z_-=(|\nabla u|^2+\sigma)^2-4(\nabla u\cdot v_*)^2.
\]
Hence
\[
 -8\ell(\ell+\sigma)\Delta b
 =2\ell Z_+^2+\sigma Z_+Z_-+8\ell\sigma Z_+,
\]
which is \eqref{eq:source}. Since $\ell,\sigma>0$ and
$Z_\pm\ge0$, the right-hand side of \eqref{eq:source}
vanishes precisely when $Z_+=0$. This is equivalent to
$v_*=-\nabla u$, or $\nabla W=-\ell\nabla u$.
\end{proof}

Although \(b\) satisfies \(\Delta b\le 0\) wherever \(\nabla W\ne0\), it need not be \(C^2\) at critical points of \(W\). We establish superharmonicity there in the viscosity sense: every \(C^2\) function \(\varphi\) touching \(b\) from below at an interior point \(P\) must satisfy \(\Delta\varphi(P)\le0\); see \cite[Definition~2.2]{CrandallIshiiLions1992}.
The proof uses the variational definition of $b$ to transfer such contact
to a suitable smooth competitor. This avoids assuming a smooth minimizing
vector at points where the quartic may have several minimizers.

\begin{lemma}\label{lem:critical}
Under the assumptions and notation of Lemma~\ref{lem:airy},
let $b$ be defined by \eqref{eq:envelope}.
If $P\in\Omega$ and $\varphi\in C^2$ near $P$ satisfies
$\varphi(P)=b(P)$ and $\varphi\le b$ near $P$, then
$\Delta\varphi(P)\le0$.
\end{lemma}

\begin{proof}
Where $\nabla W(P)\ne0$, Lemma~\ref{lem:envelope} gives $b\in C^2$
near $P$, and Proposition~\ref{prop:source} yields
$\Delta b(P)\le0$. The lower contact implies
$D^2\varphi(P)\le D^2b(P)$ as quadratic forms, and hence
$\Delta\varphi(P)\le0$.

Suppose now that $\nabla W(P)=0$. Within this proof write
$w=\nabla W$, $k_0=k(P)$, and $g_0=\nabla u(P)$.
For any $C^2$ vector field $V$ defined near $P$, put
\[
 G_V(x)=W(x)+k(x)|V(x)|^2+\frac14|V(x)|^4
                -2V(x)\cdot w(x).
\]
By \eqref{eq:envelope}, $b\le G_V$ wherever $V$ is defined. If $V(P)$ is a
minimizer of the quartic objective at $P$, then
$G_V(P)=b(P)$. In that case
\[
 \varphi\le b\le G_V,\qquad
 \varphi(P)=b(P)=G_V(P).
\]
Thus $G_V-\varphi$ has a local minimum at $P$. For every
$\xi\in\R^2$, its second derivative along $P+s\xi$ is
nonnegative at $s=0$. Hence $D^2\varphi(P)\le D^2G_V(P)$
as quadratic forms, and taking the trace gives
\begin{equation}\label{eq:competitor-test}
 \Delta\varphi(P)\le\Delta G_V(P).
\end{equation}
We choose $V$ so that the last Laplacian is nonpositive.

\smallskip
\noindent\emph{Case 1: $k_0<0$.}
Choose a constant vector $V(x)=v_0$ such that
$|v_0|^2=-2k_0$ and $v_0\cdot g_0\ge0$.
By Lemma~\ref{lem:envelope}, $v_0$ is a minimizer at $P$.
Using $\Delta W=k$, $\Delta k=-1$, and
$\Delta w=\nabla\Delta W=\nabla u$, we obtain
\[
 \Delta G_V(P)
 =k_0-|v_0|^2-2v_0\cdot g_0
 =3k_0-2v_0\cdot g_0<0.
\]
Hence \eqref{eq:competitor-test} gives $\Delta\varphi(P)<0$.

\smallskip
\noindent\emph{Case 2: $k_0=0$.}
The constant choice $V\equiv0$ minimizes the objective at $P$.
Then $G_V=W$ and
\[
 \Delta G_V(P)=k_0=0,
\]
so $\Delta\varphi(P)\le0$. This one argument covers both
$g_0=0$ and $g_0\ne0$.

\smallskip
\noindent\emph{Case 3: $k_0>0$.}
The unique minimizer at $P$ is again $v_0=0$, but we choose
the smooth competitor field $V(x)=w(x)/k_0$.
Writing $t=|w|^2$, it gives
\begin{align}
\label{eq:GV}
 G_V=W+\left(\frac{k}{k_0^2}-\frac2{k_0}\right)t
           +\frac{t^2}{4k_0^4}.
\end{align}
At $P$, $w=0$ implies $t=0$ and $\nabla t=0$.
Equation~\eqref{eq:t-laplacian} consequently gives
\[
 \Delta t(P)=k_0^2+\frac14|g_0|^4.
\]
Taking the Laplacian of \eqref{eq:GV}, all terms containing $t$ or $\nabla t$ vanish at $P$, and
\[
 \Delta G_V(P)
 =k_0-\frac{\Delta t(P)}{k_0}
 =-\frac{|g_0|^4}{4k_0}\le0.
\]
Equation~\eqref{eq:competitor-test} completes this case.

The three cases cover all possibilities when $w(P)=0$. Hence $b$ is a viscosity supersolution of $-\Delta b=0$ on $\Omega$.
\end{proof}

The strong minimum principle for viscosity supersolutions applies to \(b\) by Lemma~\ref{lem:critical}. For the boundary inequality, we use the classical Hopf barrier argument, which applies because \(b\) is \(C^2\) in an interior neighborhood of \(\partial \Omega\).

\begin{lemma}\label{lem:minimum}
Under the assumptions of Lemma~\ref{lem:airy},
let $b$ be defined by \eqref{eq:envelope} and set
$m_b:=\min_{\partial\Omega}b$.
Then $b\ge m_b$ on $\overline\Omega$,
with strict inequality in $\Omega$ unless $b$ is constant.
\end{lemma}

\begin{proof}
By Lemma~\ref{lem:critical}, $b$ is a continuous viscosity supersolution
of $-\Delta b=0$. The strong minimum principle
\cite[Proposition~4.9]{CaffarelliCabre1995} forces $b$ to be constant if
its global minimum is attained in $\Omega$. Compactness of
$\overline\Omega$ then gives both assertions.
\end{proof}

\begin{lemma}\label{lem:hopf}
Under the assumptions and notation of Lemma~\ref{lem:airy},
let $b$ be defined by \eqref{eq:envelope}.
If $b$ is nonconstant, then at every
$P\in\partial\Omega$ where
$b(P)=\min_{\partial\Omega}b$,
the function $b$ is $C^2$ up to the boundary near $P$,
and $b_\nu(P)<0$, where $\nu$ denotes the outward unit normal.
\end{lemma}

\begin{proof}
Let $P\in\partial \Omega$ satisfy $b(P)=\min_{\partial \Omega} b$, and let $u$ be the solution to \eqref{eq:system}. Since $u>0$ on $\partial \Omega$, the boundary identity
\eqref{eq:airy-boundary} gives
\[
|\nabla W|=\frac{\beta u^2}{4}>0
\qquad\text{on }\partial \Omega.
\]
By continuity and the $C^{2,\alpha}$ regularity of $\partial \Omega$, there exists a small interior tangent disk $B$ at $P$ such that $|\nabla W|$ remains nonzero.
Hence,
by Lemma~\ref{lem:envelope} and
Proposition~\ref{prop:source}, $b$ is $C^2$
and satisfies $\Delta b\le0$ in $B$. Since $b$ is nonconstant, Lemma~\ref{lem:minimum}
implies that $b>b(P)$ in $B$.  Applying the Hopf
boundary point lemma
\cite[Lemma~3.4]{GilbargTrudinger2001}
to $b$ on $B$, we obtain $b_\nu(P)<0$.
\end{proof}

\section{Boundary analysis and rigidity}\label{sec:boundary}

By Proposition~\ref{prop:constant-envelope}, a boundary
minimum of $b$ is attained at a point where $u_s=0$.
To apply Hopf's boundary point lemma, it remains to
determine the sign of $b_\nu$ at such points.
The boundary normalization \eqref{eq:airy-boundary} gives
\[
\nabla W=\frac{\beta u^2}{4}\nu,\qquad
t=|\nabla W|^2=\frac{\beta^2u^4}{16}>0
\quad\text{on }\partial\Omega.
\]
Thus the differentiable representation of $b$ in
Lemma~\ref{lem:envelope} applies near $\partial\Omega$.

\begin{lemma}\label{lem:boundary}
Under the assumptions of Lemma~\ref{lem:airy},
let $b$ be defined by \eqref{eq:envelope}, and
write $\ell=\ell(k,t)$ and $\sigma=\sigma(k,t)$,
as in Lemma~\ref{lem:envelope}. Then
\begin{equation}\label{eq:normal-b}
b_\nu
=\frac{\beta u^2}{8\ell}
\left[
\sigma\left(1-\frac{4\ell}{u}\right)^2-u_s^2
\right]
\qquad\text{on }\partial\Omega.
\end{equation}
In particular, $b_\nu\ge0$ at every boundary point
where $u_s=0$.
\end{lemma}

\begin{proof}
By \eqref{eq:t-gradient} and
$\nabla W=(\beta u^2/4)\nu$ on $\partial\Omega$,
we have
\[
t_\nu=\frac{\beta u^2}{4}
\left(k+\frac{u_s^2-\beta^2u^2}{2}\right).
\]
Using \eqref{eq:differential-b}, the Robin boundary
condition $u_\nu=-\beta u$, and
$k=\ell-\sigma/2$, we obtain
\[
\begin{aligned}
b_\nu
&=W_\nu+\sigma u_\nu-\frac{t_\nu}{\ell}\\
&=\frac{\beta u^2}{4}-\beta u\sigma
-\frac{\beta u^2}{4\ell}
\left(k+\frac{u_s^2-\beta^2u^2}{2}\right)\\
&=\frac{\beta u^2}{8\ell}
\left(\sigma-\frac{8\ell\sigma}{u}
+\beta^2u^2-u_s^2\right).
\end{aligned}
\]
Since
\[
\ell^2\sigma=t=\frac{\beta^2u^4}{16}
\qquad\text{on }\partial\Omega,
\]
this simplifies to \eqref{eq:normal-b}.
The conclusion follows from $\ell>0$ and
$\sigma>0$ on $\partial\Omega$.
\end{proof}

We now complete the proof of Theorem~\ref{thm:main}.

\begin{proof}[Proof of Theorem~\ref{thm:main}]
Suppose that \eqref{eq:system} admits a solution,
and let $b$ be defined by \eqref{eq:envelope}.
By Proposition~\ref{prop:constant-envelope},
there exists $P\in\partial\Omega$ at which
$b$ attains its boundary minimum and
$u$ attains its boundary minimum or maximum.
In particular, $u_s(P)=0$, so
Lemma~\ref{lem:boundary} gives $b_\nu(P)\ge0$.
If $b$ were nonconstant, Lemma~\ref{lem:hopf}
would yield $b_\nu(P)<0$, a contradiction.
Thus $b$ is constant on $\overline\Omega$,
and Proposition~\ref{prop:constant-envelope}
implies that $\Omega$ is a disk.
Conversely, every disk admits the radial
Robin torsion solution \eqref{eq:disk},
which satisfies \eqref{eq:system}.
\end{proof}
%==========
\begin{remark}\label{rem:profile-rigidity}
If $b$ is constant, Proposition~\ref{prop:source} gives, on
$U:=\{\nabla W\ne0\}$, the relations
$\nabla W=-\ell g$, $\sigma=|g|^2$, and $\ell=k+|g|^2/2$, where
$g=\nabla u$. Put $\mathsf A=D^2u+\Id/2$. Then
$\nabla\ell=g/2+\mathsf A g$ and
$D^2W=\ell\Id/2-g\otimes g/2$.
Differentiating $\nabla W=-\ell g$ and comparing gives
\[
 \ell\mathsf A+g\otimes(\mathsf A g)=0,
 \qquad \ell|\mathsf A|^2+|\mathsf A g|^2=0.
\]
The second identity is the Frobenius product of the first with $\mathsf A$;
since $\ell>0$, $\mathsf A=0$ on $U$. This set $U$ is dense, for
$\nabla W=0$ on an open set would give $k=\Delta W=0$ and $u=-C_\beta(\Omega)$
there, contradicting $-\Delta u=1$. Hence $D^2u=-\Id/2$ throughout
$\Omega$. The resulting quadratic torsion profile, together with the Robin boundary condition, provides an alternative route to domain rigidity.
\end{remark}

\begin{remark}
\label{rem:tensor-comparison}
Obata's conformal rigidity
argument~\cite{Obata1971} and its CR analogue due to
Jerison and Lee~\cite{JerisonLee1988} use divergence
identities with nonnegative tensorial terms to obtain
rigidity. See~\cite{LinMaOu2024} for a trace-free-Hessian
formulation and~\cite{MaWu2024,MaOuWu2025} for related
invariant-tensor constructions.
A further comparison is provided by the planar Liouville
argument in~\cite[Section~6.3]{GuiLiWangXu2026}.
For $-\Delta q=e^{2q}$, set $z=e^{-q}$ and consider
the variational stress tensor
\[
\mathsf S_q
=\frac12\bigl(|\nabla q|^2-e^{2q}\bigr)\Id
-\nabla q\otimes\nabla q.
\]
The corrected tensor
\[
\mathsf Q
=\mathsf S_q+D^2q-(\Delta q)\Id
=-\frac1z\left(
D^2z-\frac{\Delta z}{2}\Id
\right)
\]
is divergence-free and trace-free. In two dimensions,
its components are therefore harmonic. If they are
bounded on $\R^2$, then $\mathsf Q$ is constant,
yielding a linear Hessian system for $z$.

In our Robin torsion problem, the divergence-free
stress tensor $\mathsf T$ instead gives rise to a
normalized Airy stress function $W$ encoding the stationarity
condition on the boundary. A quartic envelope calibrated
on disks then produces a superharmonic auxiliary function
$b$. Its boundary structure and Hopf's lemma force $b$
to be constant, while the equality case yields the
Hessian rigidity established in
Remark~\ref{rem:profile-rigidity}.
\end{remark}

%=======
\section{Rigidity for semilinear Robin problems}\label{sec:semilinear}

In this section, we extend the method developed for Theorem~\ref{thm:main} to semilinear equations with positive, nonincreasing source terms. This extension addresses the planar case of the question raised in \cite[Section~4.1, Example~2]{BandleWagner2015}.

\begin{theorem}\label{thm:semilinear}
Let $\Omega\subset\R^2$ be a bounded simply connected $C^{2,\alpha}$
domain, $0<\alpha<1$, and $\beta>0$. Let $g\in C^1(I)$ satisfy
$g>0$ and $g'\le0$ on an open interval $I$ containing
$u(\overline\Omega)$. If a positive $u\in C^{2,\alpha}(\overline\Omega)$
solves
\begin{equation}\label{eq:semi-system}
\begin{cases}
-\Delta u=g(u)&\text{in }\Omega,\\
\partial_\nu u+\beta u=0&\text{on }\partial\Omega,\\
-\beta^2u^2+\dfrac12|\nabla u|^2
+\dfrac\beta2\kappa u^2-G(u)=C&\text{on }\partial\Omega,
\end{cases}
\end{equation}
where $G'=g$ and $C\in\R$, then $\Omega$ is a disk and $u$ is radial.
\end{theorem}

\begin{remark}
The monotonicity assumption \(g'\le0\) serves two purposes. It ensures uniqueness for the semilinear Robin problem, and hence radial symmetry of solutions on disks. It also provides the monotonicity properties of the family of radial functions required to construct superharmonic envelopes.

\end{remark}
	
\begin{proof}
The approach and structure of this proof are the same as those for Theorem~\ref{thm:main}; therefore, we will break the proof down into several steps, following the format used in the previous section.

\smallskip\noindent
\emph{\textbf{Step 1: The normalized Airy stress function.}}
As the same argument in Lemma~\ref{lem:airy}, the divergence-free tensor
$\bigl(|\nabla u|^2/2-G(u)\bigr)\Id-\nabla u\otimes\nabla u$
provides a normalized Airy stress function $W$ satisfying, with $k=G(u)+C$,
\begin{equation}\label{eq:semi-airy}
D^2W=\frac12\left[\left(k+\frac12|\nabla u|^2\right)\Id
-\nabla u\otimes\nabla u\right]\quad\text{in }\Omega; \quad
W=0,\quad\nabla W=\frac\beta4u^2\nu\ \text{on }\partial\Omega.
\end{equation}
The Hopf lemma and the Robin condition imply $u>0$ on $\overline\Omega$.

\smallskip\noindent
\emph{\textbf{Step 2: The disk-calibrated envelope.}}

Extend $g$ and its primitive $G$ outside the range
of $u$, without changing them nearby, so that
$0<m_g\le g\le M_g$ and $g'\le0$ on $\R$.
For $t\in\R$ and $R>0$, let $U_{t,R}$ be the regular radial solution of
$U''+r^{-1}U'+g(U)=0$ with $U'(0)=0$ and $U(R)=t$, and set
\[
P(t,R)=-U_{t,R}'(R),\qquad
\sigma(R)=R^2/4,\qquad h(t,\sigma)=2P(t,R)/R.
\]
The function \(h\) extends continuously differentiably to \(\sigma=0\), with \(h(t,0)=g(t)\).

To justify this definition, shooting from $U(0)=a>t$ gives
$m_{g} r/2\le -U'(r)\le M_{g} r/2$.
Writing $p=-U'$ and $z=\partial_aU$, we have
$(rz')'=-rg'(U)z$, $z(0)=1$, and hence
$z>0$, $z'\ge0$. Consequently, the hitting
radius satisfies $R_a=z(R)/p(R)>0$
and ranges from $0$ to $\infty$.
Furthermore,  the identity
$\bigl[r(p'z-pz')\bigr]'=pz/r>0$ yields $P_R>0$ (at fixed $t$).
The endpoint equation further gives
$P_R=g(t)-P/R+PP_t$. Consequently,
\begin{equation}\label{eq:semi-h-identities}
h_t=-\frac{2z'(R)}{Rz(R)}\le0,\quad
E:=\partial_\sigma(\sigma h^2)=\frac{4P P_R}{R}>0,\quad
h+\sigma h_\sigma=g(t)+\sigma hh_t=:A.
\end{equation}

To see how the envelope is determined, consider a radial solution on a disk centered at the origin, and set \(v_0=x/2\) and \(\sigma=|v_0|^2\). Then
$\nabla u=-hv_0$ and \eqref{eq:semi-airy} gives
$\nabla W=(G(u)+C+\sigma h^2/2)v_0$.
For an envelope $b=W+\min_v\{\mathcal H(u,|v|^2)-2v\cdot\nabla W\}$,
requiring $v_0$ to minimize and $\nabla b=0$ on disks suggests the
\emph{two calibration equations}
\begin{equation}\label{eq:semi-calibration}
\mathcal H_\sigma=G(t)+C+\frac12\sigma h^2,\qquad
\mathcal H_t=\sigma h.
\end{equation}
Their compatibility is exactly the last identity in
\eqref{eq:semi-h-identities}. Thus, with $\mathcal H(t,0)=0$,
\begin{equation}\label{eq:semi-envelope}
\mathcal H(t,\sigma)=\sigma(G(t)+C)
 +\frac12\int_0^\sigma\eta h(t,\eta)^2\,d\eta,
\qquad
b=W+\min_{v\in\R^2}\{\mathcal H(u,|v|^2)-2v\cdot\nabla W\}.
\end{equation}
In particular, $\mathcal H_{tt}=\sigma h_t\le0$
and $\mathcal H_{\sigma\sigma}=E/2>0$.
Uniform coercivity gives $b\in C(\overline\Omega)$.
On $\{\nabla W\ne0\}$, the minimizer $v$ is
unique and $C^1$, and $b$ is $C^2$. Setting
$\sigma=|v|^2$, $\ell=\mathcal H_\sigma(u,\sigma)>0$,
and $d=\ell+\sigma E>0$, we have
$\nabla W=\ell v$ and
\begin{equation}\label{eq:semi-gradient}
\nabla b=(\mathcal H_t+\nabla u\cdot v)\nabla u
+\frac12(\sigma h^2-|\nabla u|^2)v.
\end{equation}
Differentiation as in Proposition~\ref{prop:source} gives
\begin{equation}\label{eq:semi-laplacian}
\begin{aligned}
\Delta b=-\frac{|\nabla u+hv|^2}{8\ell d}\bigl[&
\sigma E|\nabla u-hv|^2
+2\ell|\nabla u+(2A-h)v|^2+8\ell\sigma(A^2-dh_t)\bigr]\le0.
\end{aligned}
\end{equation}
At $\nabla W=0$, the lower-test-function argument of
Lemma~\ref{lem:critical} still applies: use a constant minimizing
vector pointing along $\nabla u$ if $k<0$, the zero vector if $k=0$,
and $v(x)=\nabla W(x)/k(P)$ near the contact point $P$ if $k(P)>0$.
The respective test-function Laplacians at $P$ are negative,
zero, and $-|\nabla u(P)|^4/(4k(P))\le0$.
Thus $b$ is superharmonic in the viscosity sense on $\Omega$.

\smallskip\noindent
\emph{\textbf{Step 3: Characterization of disks.}}
By \eqref{eq:semi-airy}, on $\partial\Omega$ we have
\[
b=\mathcal B(u),\qquad
\mathcal B(t)=\min_{s\ge0}F(t,s),\qquad
F(t,s)=\mathcal H(t,s^2)-\frac\beta2t^2s.
\]
For every \(t>0\), each minimum is attained at \(s>0\), since \(F\) is coercive and \(F_s(t,0)=-\beta t^2/2<0\). Moreover,
$F_{tt}(t,s)=s^2h_t-\beta s<0$ for $s>0$.
Indeed, if $s_\theta>0$ minimizes at
$t_\theta=\theta t_1+(1-\theta)t_2$, the concavity of
$\mathcal H(\cdot,s_\theta^2)$ implies
\[
\mathcal B(t_\theta)\ge
\theta\mathcal B(t_1)+(1-\theta)\mathcal B(t_2)
+\frac\beta2s_\theta\theta(1-\theta)(t_1-t_2)^2.
\]
Thus $\mathcal{B}$ is strictly concave on $(0,\infty)$.
As in Proposition~\ref{prop:constant-envelope}, the constancy of \(b\) forces \(u\) to be constant on \(\partial\Omega\), and hence \(\kappa\) is constant. Thus \(\Omega\) is a disk.
Conversely, let $\Omega=B_R(0)$. Since \(g\) is nonincreasing and \(\beta>0\), the Robin problem has at most one solution by the standard energy argument (see \cite[Theorem~1]{MoralesRodrigoSuarez2006}). As the problem is invariant under rotations of \(B_R(0)\), \(u\circ Q\) is also a solution for every rotation \(Q\). Uniqueness therefore implies \(u=u\circ Q\), and hence \(u=U(r)\). The calibration vector $v_0=x/2$ then satisfies
$\nabla u=-hv_0$, $\nabla W=\ell(r)v_0$, and
\[
\ell(r)=G(U(r))+C+\frac12U'(r)^2,\qquad
\ell'(r)=-\frac{U'(r)^2}{r}\le0,\qquad
\ell(R)=\frac{\beta U(R)^2}{2R}>0.
\]
Since $\mathcal H_{\sigma\sigma}>0$ and $\ell>0$,
$v_0$ is the unique global minimizer in \eqref{eq:semi-envelope}.
Equation~\eqref{eq:semi-gradient} then gives $\nabla b=0$.
We have proved
\begin{equation}\label{eq:semi-equivalence}
 b\equiv\mathrm{const}\text{ on }\overline\Omega
 \quad\text{if and only if }\quad \Omega\text{ is a disk}.
\end{equation}

\smallskip\noindent
\emph{\textbf{Step 4: Rigidity.}}
On $\partial\Omega$ the minimizer has the form $v=s\nu$, $s>0$.
Equation~\eqref{eq:semi-gradient} and the Robin condition yield
\[
 b_\nu=\frac s2\bigl[(sh(u,s^2)-\beta u)^2-u_s^2\bigr].
\]
Strict concavity of $\mathcal B$ gives a boundary minimum of $b$
at a boundary extremum of $u$, where $u_s=0$ and $b_\nu\ge0$.
If $b$ were nonconstant, the strong minimum principle for
continuous viscosity supersolutions and the Hopf lemma would
instead give $b_\nu<0$ there. These principles apply as in Lemmas~\ref{lem:minimum}--\ref{lem:hopf}. Indeed, since \(|\nabla W|=\beta u^2/4>0\) on \(\partial\Omega\), continuity ensures that \(\nabla W\ne0\) in a neighborhood of \(\partial\Omega\).
Thus $b$ is constant, and \eqref{eq:semi-equivalence} finishes
the proof. 
\end{proof}
\begin{remark}
The construction reduces to the original quartic envelope \eqref{eq:envelope} when \(g\equiv1\). Indeed, choosing \(G(t)=t\), we have \(h\equiv1\) and \(\mathcal H(t,\sigma)=(t+C)\sigma+\sigma^2/4\).
\end{remark}

\section{Curvature-based rigidity via Weinberger's
\texorpdfstring{$P$}{P}-function}
\label{sec:P-function}
We conclude by revisiting Weinberger's classical
$P$-function $\mathcal P=|\nabla u|^2+u$ to illustrate
both the strengths and the limitations of this approach
for the Robin overdetermined problem.
Although $\mathcal P$ is subharmonic, the Robin condition
alone does not determine the sign of its outward normal
derivative at a boundary maximum, preventing a direct
application of Hopf's lemma. We show that the stationary
condition, together with the multiplier and mass identities,
overcomes this difficulty under $\kappa+\beta\ge0$.
This yields a complementary rigidity result without
requiring simple connectivity or any additional restriction
on $C_\beta(\Omega)$.

In contrast, the quartic-envelope method used to prove
Theorem~\ref{thm:main} obtains the required boundary sign
without a curvature restriction, but relies on the
global construction of the normalized Airy stress function.

On each boundary component, we retain the convention
$\tau=J\nu$, where $J$ is the clockwise quarter-turn
rotation and $\nu$ is the outward unit normal.
The signed curvature $\kappa$ is defined by
\eqref{eq:frenet}.

\begin{theorem}
\label{thm:curvature-rigidity}
Let $\Omega\subset\R^2$ be a bounded connected $C^{2,\alpha}$ domain, where
$0<\alpha<1$, and let $\beta>0$. Suppose that the Robin torsion function
$u_{\Omega, \beta}$ of $\Omega$ satisfies \eqref{eq:system}. If
\begin{equation}\label{eq:curvature-threshold}
 \kappa+\beta\geq0
 \qquad\text{on }\partial\Omega,
\end{equation}
then $\Omega$ is a disk.
\end{theorem}

Since $\kappa\geq0$ on the boundary of a convex domain,
Theorem~\ref{thm:curvature-rigidity} immediately yields the following
consequence.

\begin{corollary}
\label{cor:convex}
Every bounded convex $C^{2,\alpha}$ planar domain for which
\eqref{eq:system} admits a solution is a disk.
\end{corollary}

\begin{proof}[Proof of Theorem~\ref{thm:curvature-rigidity}]
Let $u$ be the solution to \eqref{eq:system}.
Consider the $\mathcal{P}$-function defined as
\begin{equation*}
 \mathcal{P}=|\nabla u|^2+u.
\end{equation*}
Since $\Delta u=-1$, direct differentiation yields
\begin{equation}\label{eq:weinberger-laplacian}
 \Delta \mathcal{P}=2|D^2u|^2-1
 =2\left|D^2u+\frac12\Id\right|^2\ge0.
\end{equation}
Suppose that $\mathcal{P}$ is not constant. By the strong maximum principle,
$\mathcal{P}$ attains its maximum over $\overline\Omega$ at some
$p\in\partial \Omega$. Since $\mathcal{P}\in C^2(\Omega)\cap C^1(\overline\Omega)$
and the $C^{2,\alpha}$ boundary satisfies the interior sphere
condition, the Hopf boundary point lemma implies
\begin{equation}\label{eq:weinberger-Hopf}
 \mathcal{P}_\nu(p)>0.
\end{equation}

On $\partial \Omega$, the Robin boundary condition gives
\begin{equation}\label{eq:weinberger-boundary}
 \nabla u=u_s\tau-\beta u\nu,
 \qquad \mathcal{P}=u_s^2+\beta^2u^2+u.
\end{equation}
For vectors $\xi,\eta$, write $u_{\xi\eta}=D^2u(\xi,\eta)$ for a Hessian
component. In contrast, $u_{ss}=\partial_s^2(u|_{\partial\Omega})$ is
the second derivative of the boundary trace along arclength.
Differentiating the boundary trace and the Robin condition
tangentially, and using \eqref{eq:frenet}, we obtain
\begin{align}
 u_{\tau\tau}&=u_{ss}-\beta\kappa u,
 \label{eq:weinberger-tt}\\
 u_{\tau\nu}&=-(\kappa+\beta)u_s,
 \label{eq:weinberger-tn}\\
 u_{\nu\nu}&=-1-u_{ss}+\beta\kappa u.
 \label{eq:weinberger-nn}
\end{align}
Here the last identity follows from $\Delta u=-1$.

At the maximum point $p\in \partial \Omega$ of $\mathcal{P}$, we claim that $u_s(p)\ne0$. Suppose instead that $u_s(p)=0$.
Let $M=\max_{\partial \Omega} u$ and choose $p_M\in\partial \Omega$ such that
$u(p_M)=M$. Since $u_s(p_M)=0$, maximality of $\mathcal{P}(p)$ gives
\[
 \beta^2u(p)^2+u(p)=\mathcal{P}(p)\ge \mathcal{P}(p_M)
 =\beta^2M^2+M.
\]
The function $t\mapsto\beta^2t^2+t$ is strictly increasing on
$(0,\infty)$; as $0<u(p)\le M$, it follows that $u(p)=M$.
In particular, $u_{ss}(p)\le0$.

The identities \eqref{eq:mass-energy} and \eqref{eq:constant} give
\begin{equation*}
 -C_{\beta}(\Omega)=\frac{1}{2|\Omega|}\int_\Omega \mathcal{P}\dx
   +\frac{\beta}{4|\Omega|}\int_{\partial \Omega} u^2\ds.
\end{equation*}
Since $\mathcal{P}\le \mathcal{P}(p)=\beta^2M^2+M$ in $\Omega$, and
\[
 \beta\int_{\partial \Omega} u^2\ds
 \le\beta M\int_{\partial \Omega} u\ds=M|\Omega|,
\]
we deduce that
\begin{equation}\label{eq:weinberger-multiplier-bound}
 -C_{\beta}(\Omega)\le\frac12(\beta^2M^2+M)+\frac14M
 =\frac12\beta^2M^2+\frac34M.
\end{equation}
On the other hand, the stationary boundary identity
\eqref{eq:scalar} at $p$, where $u_s(p)=0$ and $u(p)=M$, yields
\[
 -C_{\beta}(\Omega)=\frac12\beta^2M^2-\frac\beta2\kappa(p)M^2+M.
\]
Comparing with \eqref{eq:weinberger-multiplier-bound} shows that
\begin{equation}\label{eq:weinberger-pinching}
 \beta\kappa(p)M\ge\frac12.
\end{equation}
From \eqref{eq:weinberger-pinching}, \eqref{eq:weinberger-nn} and $u_{ss}(p)\le0$,
\[
 u_{\nu\nu}(p)+\frac12
 =-\frac12-u_{ss}(p)+\beta\kappa(p)M\ge0.
\]
Because $u_s(p)=0$ and $u_\nu(p)=-\beta M<0$, we obtain
\[
 \mathcal{P}_\nu(p)
 =2u_\nu(p)\left(u_{\nu\nu}(p)+\frac12\right)
 \le0,
\]
contrary to \eqref{eq:weinberger-Hopf}. Hence $u_s(p)\ne0$.

Differentiating \eqref{eq:weinberger-boundary} along $\partial \Omega$,
we have
\[
 0=\mathcal{P}_s(p)
 =u_s(p)\bigl(2u_{ss}(p)+2\beta^2u(p)+1\bigr).
\]
Consequently,
\[
 u_{ss}(p)=-\beta^2u(p)-\frac12,
 \qquad
 u_{\nu\nu}(p)+\frac12
 =\beta(\kappa(p)+\beta)u(p),
\]
where the second identity follows from \eqref{eq:weinberger-nn}.
Using \eqref{eq:weinberger-boundary} and
\eqref{eq:weinberger-tn}, we therefore find
\begin{align*}
 \mathcal{P}_\nu(p)
 &=2D^2u(\nabla u,\nu)(p)+u_\nu(p)\\
 &=2u_s(p)D^2u(\tau,\nu)(p)
   -2\beta u(p)\left(u_{\nu\nu}(p)+\frac12\right)\\
 &=-2(\kappa(p)+\beta)
     \bigl(u_s(p)^2+\beta^2u(p)^2\bigr)\le0,
\end{align*}
where the last inequality uses \eqref{eq:curvature-threshold}.
This again contradicts \eqref{eq:weinberger-Hopf}.
We conclude that $\mathcal{P}$ is constant in $\Omega$.

By \eqref{eq:weinberger-laplacian}, this implies
$D^2u=-\Id/2$ throughout $\Omega$. Since $\Omega$ is connected,
there exist $c\in\R$ and $z\in\R^2$ such that
\begin{equation}\label{eq:weinberger-quadratic}
 u(x)=c-\frac14|x-z|^2 \qquad (x\in\overline\Omega).
\end{equation}
The Robin condition yields
\begin{equation}\label{eq:weinberger-radial-robin}
 (x-z)\cdot\nu=2\beta u>0\qquad\text{on }\partial \Omega.
\end{equation}
Let $r_-:=\min_{x\in\partial\Omega}|x-z|$ and
$r_+:=\max_{x\in\partial\Omega}|x-z|$, attained at boundary points
$x_-$ and $x_+$, respectively.
By \eqref{eq:weinberger-radial-robin}, $x_\pm\ne z$.
At these extrema, tangential differentiation of $|x-z|^2$
shows that $x_\pm-z$ is normal to $\partial \Omega$.
Its positive outward normal component in
\eqref{eq:weinberger-radial-robin} therefore gives
$x_\pm-z=r_\pm\nu(x_\pm)$.
Combining \eqref{eq:weinberger-quadratic} and
\eqref{eq:weinberger-radial-robin}, we obtain
\[
 r_\pm+\frac\beta2r_\pm^2=2\beta c.
\]
As the function $r\mapsto r+\beta r^2/2$ is strictly increasing
for $r\ge0$, we have $r_-=r_+=:R>0$.
Thus $\partial \Omega\subset\partial B_R(z)$. Since $\Omega$ is a
nonempty bounded domain and its boundary is contained in
$\partial B_R(z)$, necessarily $\Omega=B_R(z)$.
\end{proof}

\begin{remark}
The argument in the proof of Theorem~\ref{thm:curvature-rigidity} is specific to dimension two. Indeed, 
the planar boundary has only one tangential direction, so once \(u_s(p)\ne0\) has been established at a boundary maximum of the nonconstant \(\mathcal P\)-function, the condition \(\mathcal P_s(p)=0\) determines \(u_{ss}(p)\). In higher dimensions, additional transverse tangential
components remain uncontrolled.
\end{remark}

	\section*{Statement on the use of AI tools}
	During the preparation of this manuscript, the authors used ChatGPT (OpenAI) for language and presentation assistance and as an interactive tool for exploratory mathematical discussion. These interactions served to explore and refine possible lines of argument. All mathematical statements, proofs, and computations in the final manuscript were independently verified by the authors, who take full responsibility for the article.

\end{document}